\documentclass[AMA,STIX1COL]{WileyNJD-v2}
\usepackage{moreverb}
\usepackage{subfigure} 
\usepackage{graphicx}

\newcommand\BibTeX{{\rmfamily B\kern-.05em \textsc{i\kern-.025em b}\kern-.08em
T\kern-.1667em\lower.7ex\hbox{E}\kern-.125emX}}

\articletype{Article Type}%

\received{<day> <Month>, <year>}
\revised{<day> <Month>, <year>}
\accepted{<day> <Month>, <year>}

\begin{document}


\title{Locally Adaptive Non-Hydrostatic Shallow Water Extension for Moving Bottom-Generated Waves}

\author[1,2]{Kemal Firdaus*}

\author[1,2]{J\"orn Behrens}

\authormark{FIRDAUS and BEHRENS}

\address[1]{Department of Mathematics, Universit\"at Hamburg, Bundesstrasse 53-55, 20146 Hamburg, Germany}

\address[2]{Center for Earth System Research and Sustainability (CEN), Universit\"at Hamburg, Bundesstrasse 53-55, 20146 Hamburg, Germany}

\corres{*Kemal Firdaus, Department of Mathematics, Universit\"at Hamburg, Bundesstrasse 53-55, 20146 Hamburg, Germany.\\ \email{kemal.firdaus@uni-hamburg.de}}


\abstract[Abstract]{We propose a locally adaptive non-hydrostatic model and apply it to wave propagation generated by a moving bottom. This model is based on the non-hydrostatic extension of the shallow water equations (SWE) with a quadratic pressure relation, which is suitable for weakly dispersive waves. The approximation is mathematically equivalent to the Green-Naghdi equations. Applied globally, the extension requires solving an elliptic system of equations in the whole domain at each time step. Therefore, we develop an adaptive model that reduces the application area of the extension and by that the computational time. The elliptic problem is only solved in the area where the dispersive effect might play a crucial role. To define the non-hydrostatic area, we investigate several potential criteria based on the hydrostatic SWE solution. We validate and illustrate how our adaptive model works by first applying it to simulate a simple propagating solitary wave, where exact solutions are known. Following that, we demonstrate the accuracy and efficiency of our approach in more complicated cases involving moving bottom-generated waves, where measured laboratory data serve as reference solutions. 
The adaptive model yields similar accuracy as the global application of the non-hydrostatic extension while reducing the computational time by more than $50\%$.}


\maketitle

\section{Introduction}\label{sec1}
Non-hydrostatic pressure is proven to be crucial in some geophysical flow phenomena, including tsunamis caused by landslides \cite{glimsdal2013}. Such phenomena can be described as moving bottom-generated waves. The non-hydrostatic pressure yields a wave dispersion effect. Previously proposed models including wave dispersion often involve higher-order derivatives or solving elliptic equations, which demand for more computational resources compared to non-dispersive models. This motivates us to develop an approach, where the dispersive effect is locally and adaptively applied to moving bottom-generated waves.

Depth-averaged models are commonly used for large-scale ocean modeling due to their simplicity and computational efficiency. Boussinesq-type equations are a popular choice for depth-averaged models that can handle the dispersive effect. This model type has been widely used in tsunami research, especially for moving-bottom generated waves \cite{lynett2002, ashtiani2007, FUHRMAN2009747, dutykh2013, FANG2020101977}. It is derived based on asymptotic expansion on the velocity potential, where different expansion orders lead to various Boussinesq-type equations \cite{peregrine1967, nwogu1993, madsen1998, madsen2003}. However, such models involve higher-order and sometimes mixed time-space derivatives, which often cause numerical and computational challenges.

Another widely used modeling approach is given by the shallow water equations (SWE), involving only first-order derivatives, simplifying the solution. Despite their capability to describe a wide range of wave dynamics, SWE are limited to the hydrostatic assumption, excluding dispersive effects. Previous research proposed to extend the SWE by including a non-hydrostatic pressure term. To derive such an extension, approximations to the non-hydrostatic pressure and the vertical velocity profile need to be defined. Most studies assumed a linear vertical pressure relation \cite{walters2005, stelling2003, yamazaki2009} while also assuming a linear vertical velocity. Alternatively, a quadratic pressure relation was proposed, which is consistent with the linear vertical velocity assumption \cite{Jeschke2017}. This latter relation
was proven to be equivalent to Boussinesq-type equations, namely the Green-Naghdi equations \cite{Green_Naghdi_1976}, making it suitable for modeling weakly dispersive waves. Recently, this work was improved, making it solvable with a projection-based method and applicable for moving bottom-generated waves \cite{firdaus2025}.

Applying the extension globally to the SWE requires the solution of an elliptic system of equations over the whole domain at each time step, adding substantial computational effort to the non-hydrostatic SWE solution. Therefore, we propose an adaptive approach, which aims to apply the non-hydrostatic extension -- and thus solving the elliptic problem -- only locally in the particular area where the non-hydrostatic pressure plays a significant role. 

The projection-based method allows us to use the solution from the hydrostatic SWE solver as a predictor, which will then be corrected by solving the extension-related terms. With this principle in mind, we can apply the correction locally. To define these corrected regions, we investigate several potential criteria based on the hydrostatic SWE solver solutions. To validate and illustrate how our model works, we first compare it with an analytical solution of a simple solitary wave. Subsequently, we apply it to simulate moving bottom-generated waves, which are based on two previous laboratory setups. These experiments involved both vertical and horizontal movement.


\section{Mathematical Model}\label{sec2}
In this study, we use the one-dimensional non-hydrostatic SWE extension, which is also applicable for moving bottom-generated waves \cite{firdaus2025}, written as:
\begin{equation}\label{eq:swe_mass}
   h_t+(hu)_x=0,
\end{equation}
\begin{equation}\label{eq:swe_momentum_1}
    (hu)_t+\left(hu^2+\frac{ g}{2}h^2\right)_x=ghd_x-\frac{1}{\rho}(hp^{nh})_x+\frac{1}{\rho}P^{nh}|_{z=-d}d_x,
\end{equation}
\begin{equation}\label{eq:swe_momentum_2}
    (hw)_t+(huw)_x=\frac{1}{\rho}P^{nh}|_{z=-d},
\end{equation}
\begin{equation}\label{eq:swe_constraint}
    2hw +hu(2d-h)_x+2hd_t = -h(hu)_x
\end{equation}
where Eq. \ref{eq:swe_mass} is the mass conservation equation, while Eq. \ref{eq:swe_momentum_1} and \ref{eq:swe_momentum_2} are the horizontal and vertical momentum balance equations respectively. These equations are also completed with the divergence constraint Eq. \ref{eq:swe_constraint}, which is derived from a linear vertical velocity assumption. The fluid thickness is denoted as $h = \eta+d$, considering the surface elevation $\eta$ and fluid depth $d$. Meanwhile, $u$ is the depth-averaged horizontal velocity. In contrast to the hydrostatic SWE, this model also involves the depth-averaged vertical velocity $w$ and non-hydrostatic pressure $p^{nh}$. To enclose this set of equations, a pressure relation needs to be defined. This study uses a quadratic relation that was further developed for moving bottom cases \cite{firdaus2025}, written as:
\begin{equation}\label{eq:pressure_rev}
    P_{-d}^{nh}=\frac{6}{4+d_x^2}p^{nh}+\frac{d_x}{4+d_x^2}(hp^{nh})_x+\phi,
\end{equation}
with $\phi = \frac{\rho h}{4+d_x^2}\left(gd_x\eta_x-d_{tt}-2ud_{xt}-u^2d_{xx}\right)$. These equations allow us to solve it with a projection method, which further could be used to obtain an adaptive model.

\section{Numerical Method}\label{sec3}
We will briefly discuss how to solve the model as described in the previous section. As mentioned before, we are interested in solving it using a projection method, also known as a predictor-corrector method. This method lets us use the hydrostatic SWE solver solution as the predictor and make a correction by solving the extension terms. The correction itself is not solved explicitly but rather by solving a system of elliptic equations, which will be elaborated on later.

\subsection{Predictor Step}
The hydrostatic SWE (Eq. \ref{eq:swe_mass}, \ref{eq:swe_momentum_1}, and \ref{eq:swe_momentum_2} without non-hydrostatic pressure terms) are discretized with a second-order Runge-Kutta discontinuous Galerkin (RKDG) method. We consider a computational domain $\Omega = [x_l,x_r]$, which is then divided into $N$ uniform intervals $I_i = [x_{i-1},x_i]$, where $i = 1,2,..,N$, assuming $x_0 = x_l$ and $x_N = x_r$. To formulate the weak DG formulation of the hydrostatic SWE, we multiply it with a test function $\phi$ and integrate it over an interval, which leads to
\begin{equation}\label{eq:DG_weak}
    \int_{I_i}\phi\boldsymbol{Q}_tdx-\int_{I_i}\phi_x\boldsymbol{F}(\boldsymbol{Q})dx+\left[\phi\boldsymbol{F}^*(\boldsymbol{Q})\right]_{x_i}^{x_{i+1}}=\int_{I_i}\phi\boldsymbol{S}(\boldsymbol{q})dx,
\end{equation}
where $\boldsymbol{Q} = (h, hu, hw)^T$ is the unknowns, while the fluxes and sources are denoted by $\boldsymbol{F}(\boldsymbol{Q}) = \left(hu,hu^2+\frac{g}{2}h^2,huw\right)^T$ and $\boldsymbol{S}(\boldsymbol{Q}) = (0, ghd_x, 0)^T$ respectively. The communication between neighboring intervals depends on the definition of the numerical interface flux $\boldsymbol{F}^*(\boldsymbol{Q})$, which is unique and based on information from both elements. A Riemann solver (in this case, the Rusanov solver) is used as an approximation.

We discretize the solution in space with a piecewise polynomial represented by the nodal Lagrange basis function, which is also used as the test function. In other word, our solution can be approximated as $\boldsymbol{Q}_h(x,t) = \sum_i(\tilde{\boldsymbol{Q}}_h)_i(t)\phi_i(x)$. The fluxes and the sources are approximated similarly. This spatial discretization results in a system of ordinary differential equations for the local degrees of freedom, such that
\begin{equation}
    \frac{d\tilde{\boldsymbol{Q}}_h}{dt} = \boldsymbol{M}^{-1}\left(\boldsymbol{D}^T\boldsymbol{F}(\tilde{\boldsymbol{Q}}_h)-\boldsymbol{M}\boldsymbol{S}(\tilde{\boldsymbol{Q}}_h)\right)-\left[\phi\boldsymbol{F}^*(\tilde{\boldsymbol{Q}}_h)\right]_{x_i}^{x_{i+1}} =: \mathcal{H}(\tilde{\boldsymbol{Q}}_h),
\end{equation}
where $\boldsymbol{M}_{ij} = \int_{I_i}\phi_i\phi_jdx$ and $\boldsymbol{D}_{ij} = \int_{I_i}\phi_i(\phi_j)_xdx$, and $\mathcal{H}$ represents the evolution operator of the SWE. This last equation can be solved with any time integrator, which in this work is a second-order Heun's scheme, written as
\begin{equation}
    \tilde{\boldsymbol{Q}}_h^{n+1} = \tilde{\boldsymbol{Q}}_h^{n} + \frac{\Delta t}{2}\left(\mathcal{H}(\tilde{\boldsymbol{Q}}_h^n) + \mathcal{H}(\tilde{\boldsymbol{Q}}_h^{\star})\right),
\end{equation} 
where the superscript $n$ stands for a discrete time step, and $\tilde{\boldsymbol{Q}}_h^{\star}$ is derived by one Euler forward step 
\begin{equation}
    \tilde{\boldsymbol{Q}}_h^{\star} = \tilde{\boldsymbol{Q}}_h^{n} + \Delta t\left(\mathcal{H}(\tilde{\boldsymbol{Q}}_h^n)\right).
\end{equation}
We mark the results from this middle step with a tilde, written as $\tilde{\boldsymbol{Q}}_h^{n+1}=(\tilde{h}_h^{n+1},(\tilde{hu})_h^{n+1},(\tilde{hw})_h^{n+1})^T$.

\subsection{Correction Step}
First of all, note that there is no non-hydrostatic term in the mass equation. Therefore, no correction is made to it, i.e. $h^{n+1} = \tilde{h}^{n+1}$. The correction is then made to the momentum equations, which are solved based on the backward Euler method, following:
\begin{equation}\label{eq:correction1}
    \frac{(hu)^{n+1}-(\tilde{hu})^{n+1}}{\Delta t} = -\frac{1}{\rho}(\tilde{h}p^{nh})_x^{n+1}+\frac{d_x^{n+1}}{\rho}\left(\frac{6}{4+(d_x^{n+1})^2}(p^{nh})^{n+1}+\frac{d_x^{n+1}}{4+(d_x^{n+1})^2}(\tilde{h}p^{nh})_x^{n+1}+\phi^{n+1}\right),
\end{equation}
\begin{equation}\label{eq:correction2}
\frac{(hw)^{n+1}-(\tilde{hw})^{n+1}}{\Delta t} = \frac{1}{\rho}\left(\frac{6}{4+(d_x^{n+1})^2}(p^{nh})^{n+1}+\frac{d_x^{n+1}}{4+(d_x^{n+1})^2}(\tilde{h}p^{nh})_x^{n+1}+\phi^{n+1}\right),
\end{equation}
with $\phi^{n+1} := \phi((d,\tilde{h},\tilde{hu})^{n+1})$. Through some manipulation, these last two equations yield an elliptic system of equations 
\begin{equation}\label{eq:elliptic1}
    (p^{nh})_x^{n+1}+\frac{2\tilde{h}_x^{n+1}-3d_x^{n+1}}{2\tilde{h}^{n+1}}(p^{nh})^{n+1}+\frac{\rho(4+(d_x^{n+1})^2)}{4\Delta t\tilde{h}^{n+1}}(hu)^{n+1} = \frac{4+(d_x^{n+1})^2}{4\tilde{h}^{n+1}}\left(\phi^{n+1} d_x^{n+1}+\frac{\rho}{\Delta t}(\tilde{hu})^{n+1}\right).
\end{equation}
\begin{equation}\label{eq:elliptic2}
    (hu)_x^{n+1}+\frac{3d_x^{n+1}-2\tilde{h}_x^{n+1}}{2\tilde{h}^{n+1}}(hu)^{n+1}+\frac{3\Delta t}{\rho \tilde{h}^{n+1}}(p^{nh})^{n+1}=-2d_t^{n+1}-\frac{2(\tilde{hw})^{n+1}}{\tilde{h}^{n+1}}-\frac{d_x^{n+1}(\tilde{hu})^{n+1}}{2\tilde{h}^{n+1}}-\frac{\Delta t(4+(d_x^{n+1})^2)}{2\rho\tilde{h}^{n+1}}\phi^{n+1},
\end{equation}
where $\phi^{n+1} := \phi((d,\tilde{h},\tilde{hu})^{n+1})$. These last two equations form an elliptic system of equations, which can be solved with the local discontinuous Galerkin (LDG) method. This method was initially developed by \cite{cockburn1998} for a time-dependent, non-linear advection-diffusion system, which was further developed to solve the Poisson equation \cite{cockburn2003, CASTILLO20061307} by rewriting it into a system of first-order equations, similar to what we have in Eq. \ref{eq:elliptic1} and \ref{eq:elliptic2}. Furthermore, the correction to the vertical momentum is done based on the obtained corrected pressure.

Let us write Eq. \ref{eq:elliptic1} and \ref{eq:elliptic2} in a shorter form and neglect the superscript for simplicity as:
\begin{equation}\label{eq:short_elliptic1}
    (p^{nh})_x+s_{11}(p^{nh})+s_{12}(hu) = f_1,
\end{equation}
\begin{equation}\label{eq:short_elliptic2}
    (hu)_x+s_{21}(hu)+s_{22}(p^{nh})=f_2.
\end{equation}
Note that $s_{11}+s_{21} = 0$ and $s_{12}>0$ in any case. These properties are important for the flux derivation, as stated in \cite{jeschke2018}. This flux is the generalized form of the previously derived flux for the Poisson equation \cite{CASTILLO20061307, cockburn2003}. The generalized flux can be expressed as:
\begin{equation}\label{eq:elliptic_flux1}
    (hu)_h^* = \{(hu)_h\}+c_{11}[p_h^{nh}]+c_{12}[u_h],
\end{equation}
\begin{equation}\label{eq:elliptic_flux2}
    {p_h^{nh}}^* = \{p_h^{nh}\}-c_{12}[p_h^{nh}]+c_{22}[u_h],
\end{equation}
where $\{q\} = (q^{-}+q^{+})/2$ and $[q] = q^{-}-q^{+}$ denote the neighboring element average and jump, with superscript $-$ and $+$ represent the value on the left and right element respectively. However, instead of setting $c_{11} = 1$ as in \cite{jeschke2018}, we set $c_{11} = \frac{1}{2}$ to achieve the "flip-flop" behavior, which characterizes the LDG method. This lets us achieve a more efficient computation as it has more zero entries in its matrix form. To handle the boundary, we use the derived zero Dirichlet boundary conditions.

\section{Preliminary Results on Global Model}\label{sec4}
Before developing the adaptive model, we investigate the validity of our numerical scheme by applying the extension globally on the whole computational domain. We begin with a simple propagating solitary wave. Thereafter, we test more complicated cases, where the waves are generated by a moving bottom, including those by vertical and horizontal movement. To quantitatively measure the accuracy of our model, we calculate the root mean square error (RMSE) of the water height to measure the error on the amplitude. Along with that, we also compute the Pearson correlation, denoted with $r$, to capture the waveform similarity. The latter index is crucial to show that our model produces the expected dispersive behavior. 

\subsection{Solitary wave}
This first benchmark considers a solitary wave propagating to the positive $x$-axis direction. Since our model is proven to be equivalent to the Serre equations for the flat bottom case, we can compare it to its analytical solution as given in \cite{seabrasantos1987}, expressed as:
\begin{equation}\label{eq:analytical}
    \eta(x,t) = a\cosh^{-2}(K(x-ct-x_0)), \quad u(x,t) = \frac{c\eta(x,t)}{d+\eta(x,t)}.
\end{equation}
This solution is suitable for a solitary wave with an amplitude of $a = 2~m$ over a constant depth $d = 10~m$, with a propagation velocity $c = \sqrt{g(d+a)}$ and a scaling factor $K = \sqrt{\frac{3a}{4d^2(d+a)}}$. In this case, the initial displacement is located at $x_0 = L/4$, where we consider a domain with a length of $L = 800~m$.

For our simulation, we simply substitute $t = 0$ to Eq. \ref{eq:analytical} as our initial conditions for the water height and the horizontal momentum. The initial vertical momentum could further be calculated from Eq. \ref{eq:swe_constraint}. We run our simulation until $T = 30~s$, where we discretize it with a time-step $\Delta t = 0.1~s$ and element width $\Delta x = 4~m$. 

\begin{figure}[h]
    \centering
    \includegraphics[width=1\linewidth]{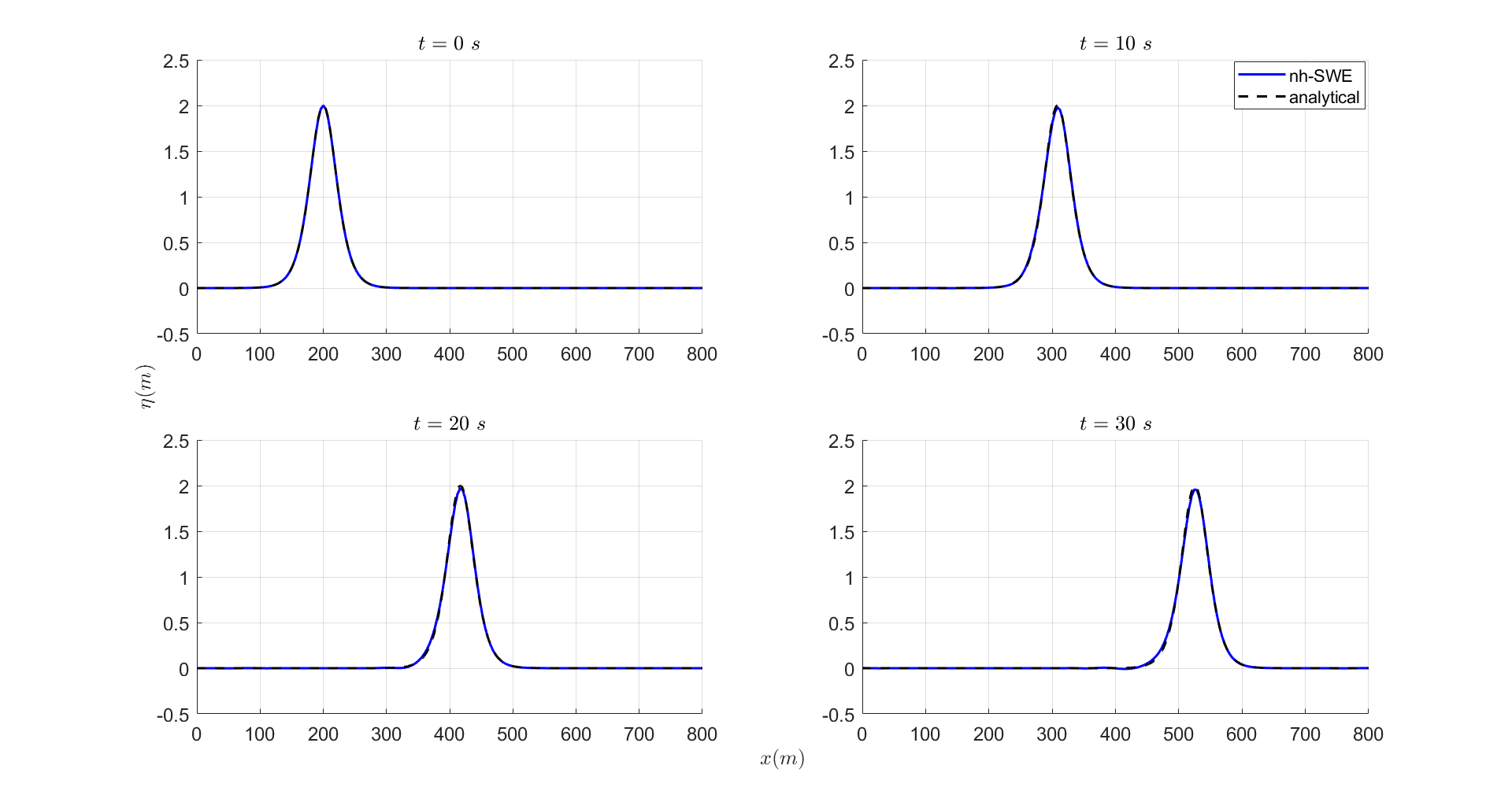}
    \caption{Comparison of our numerical simulation results (blue solid line) and its analytical solution (black dashed line) recorded at $t = 0, 10, 20, 30~s$.}
    \label{fig:global_solitary}
\end{figure}

Fig. \ref{fig:global_solitary} compares our simulation with the analytical solution at the initial time $t = 0~s$ and three timestamps $t = 10, 20, 30~s$. Our numerical results match the analytical solution very well, with an RMSE of $0.0108$ and correlation $r = 0.9997$ by the end of the simulation time. We note that the accuracy stems from the use of a quadratic pressure profile in the vertical, whereas a linear profile yields a notable inaccuracy in the amplitude for the same simulation time, as witnessed in \cite{walters2005, stelling2003, Jeschke2017}.

\subsection{Impulsive vertical thrust}
We investigate a more complicated test case, which involves a moving bottom. This benchmark is based on a laboratory experiment in \cite{Hammack_1973}, where still water is disrupted by a vertical flip movement, generating a wave. We define a wall boundary condition next to the moving plate and an absorbing boundary on the other side. Mathematically, the bottom movement can be expressed as:
\begin{equation}
    d(x,t) = h_0-\zeta_0(1-e^{-\alpha t})\mathcal{H}(b^2-x^2),
\end{equation}
where $h_0 = 0.05~m$ represents the still water depth, $|\zeta_0| = 0.005m$ is the plate maximum displacement and $b = 0.61~m$ is the plate's width. The movement of the plate follows the exponential equation, which depends on $\alpha = \frac{1.11}{t_c}$, where $t_c$ being the characteristic time that fulfills $\frac{t_c\sqrt{gh_0}}{b}=0.148$ and $0.093$ for upward and downward case respectively.

\begin{figure}[h!]
    \centering
    \includegraphics[width=1\linewidth]{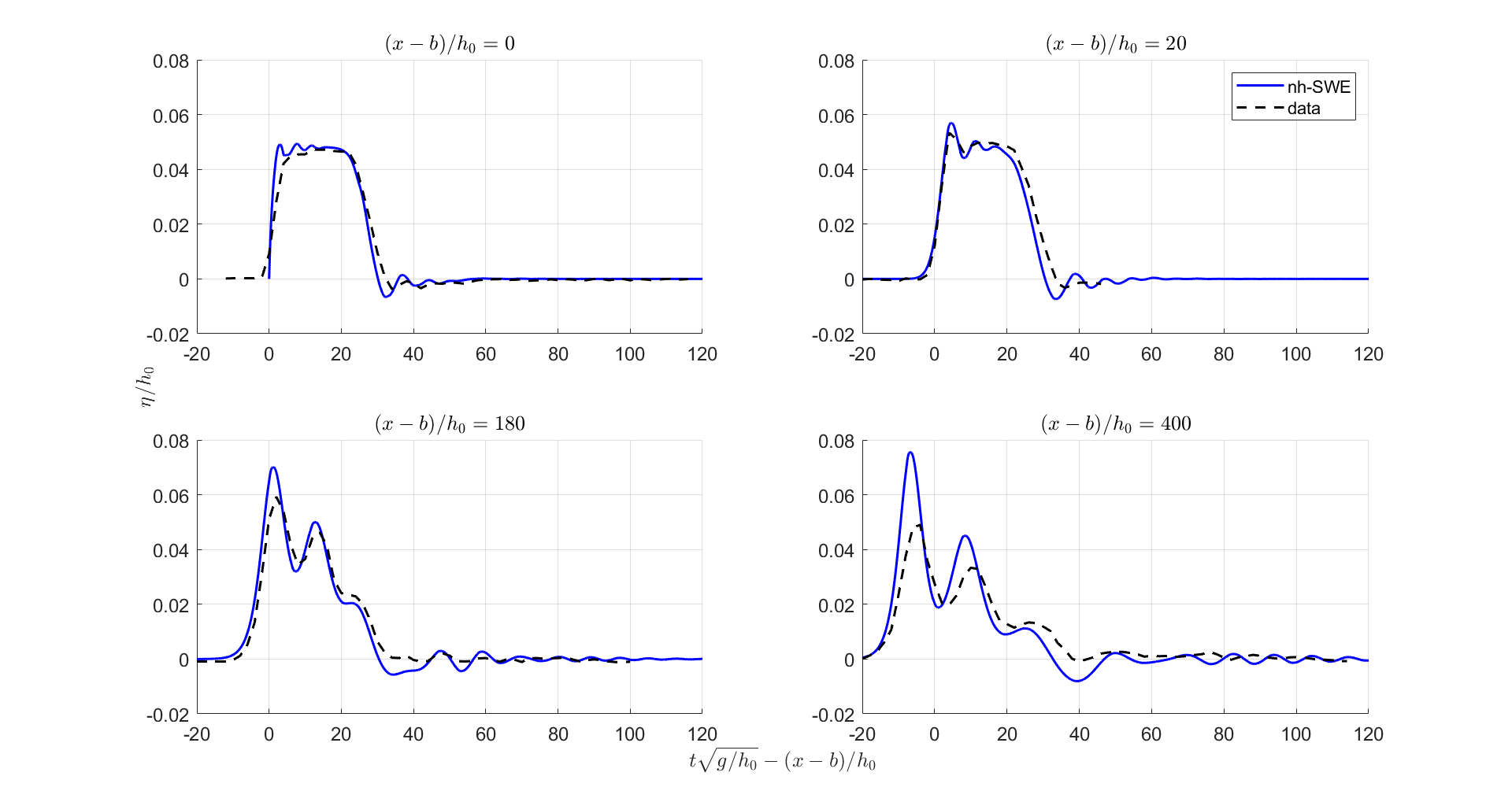}
    \caption{Comparison of the numerical simulations (blue solid line) and measured laboratory data (black dashed line) recorded at $x = 0.61, 1.61, 9.61, 20.61~m$ for the uplift case.}
    \label{fig:global_up}
\end{figure}

\begin{figure}[h!]
    \centering
    \includegraphics[width=1\linewidth]{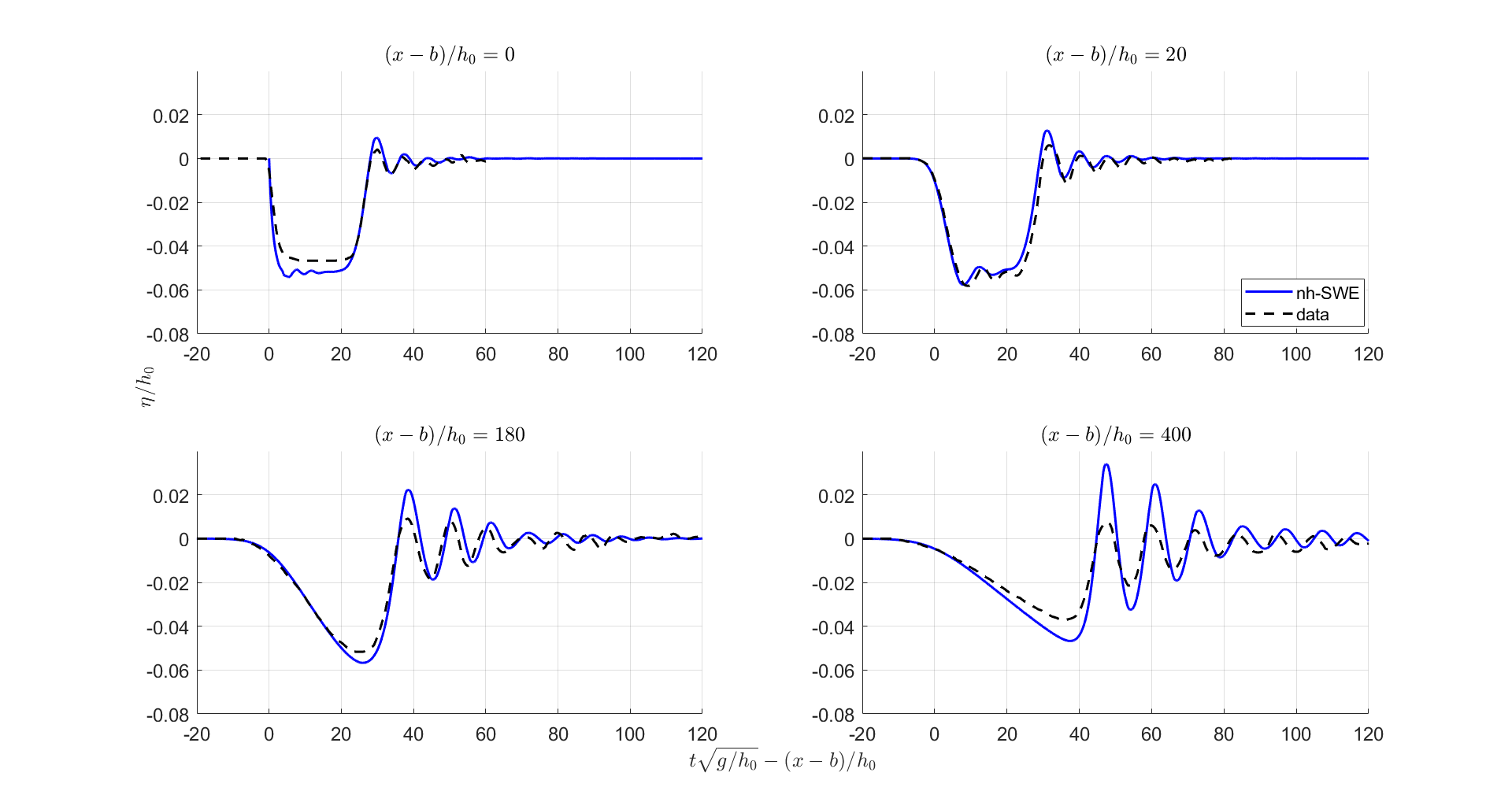}
    \caption{Comparison of the numerical simulations (blue solid line) and measured laboratory data (black dashed line) recorded at $x = 0.61, 1.61, 9.61, 20.61~m$ for the downdraft case.}
    \label{fig:global_down}
\end{figure}

To achieve a fair resolution with the measured data, we use $\Delta t = 0.01~s$ and $\Delta x = 0.025~m$. The comparison is made with the measured data at four measurement points $x = \{0.61, 1.61, 9.61, 20.61\}~m$, equivalent as $\frac{x-b}{h_0}={0, 20, 180, 400}$. The comparisons are given in Fig. \ref{fig:global_up} and \ref{fig:global_down}. All in all, it can be observed that our simulation results are in good agreement with the measured data, especially near the generation area $x=\{0.61, 1.61\}~m$. In particular, the phases of the dispersive wave trains are well captured by the simulation. Table \ref{tab:global_hammack} shows the quantitative comparisons with the measured data. From both cases, we can observe the dispersive behavior from our model, which also applies to the measured data. However, there is unavoidably a mismatch in the area away from the further measurement points $x=\{0.61, 1.61\}~m$, which is also the case in the previous studies with higher-order models \cite{FUHRMAN2009747, Hammack_1973, xin024}. This might occur due to the viscous energy losses and boundary stresses in the experiments as suggested in \cite{Hammack_1973}, which is neglected here and in all the mentioned studies. 

\begin{table}[h!]
    \centering
    \resizebox{\textwidth}{!}{%
    \begin{tabular}{|c|c|c|c|c|c|c|c|c|}
        \hline
         \multirow{2}{*}{Movement} & \multicolumn{4}{c|}{RMSE} & \multicolumn{4}{c|}{$r$}\\
         \cline{2-9}
         & $\frac{x-b}{h_0} = 0$ & $\frac{x-b}{h_0} = 20$ & $\frac{x-b}{h_0} = 180$ & $\frac{x-b}{h_0} = 400$ & $\frac{x-b}{h_0} = 0$ & $\frac{x-b}{h_0} = 20$ & $\frac{x-b}{h_0} = 180$ & $\frac{x-b}{h_0} = 400$\\
         \hline
         Uplift      & $0.0018$ & $0.0022$ & $0.0028$ & $0.0063$ 
                     & $0.9958$ & $0.9936$ & $0.9894$ & $0.9589$\\
         Downdraft   & $0.0031$ & $0.0029$ & $0.0037$ & $0.0064$
                     & $0.9955$ & $0.9902$ & $0.9829$ & $0.9609$\\

         \hline
    \end{tabular}%
    }
    \caption{Root mean square error and correlation with respect to the measured data for each case.}
    \label{tab:global_hammack}
\end{table}



\subsection{Sliding bump over flat bottom}
This case also involves a moving bottom, yet with a horizontal movement. It follows an experimental setup in \cite{whittaker2015}, exposed to a $h_0 = 0.175~m$ still water level. The wave was generated by a sliding semi-elliptical aluminum block installed in the middle of the flume. A formulation of the bottom movement was proposed \cite{Jing2020}, written as: 
\begin{equation}
    d(x,t) = h_0-H_s\left(1-\left(\frac{2(x-S(t))}{L_s}\right)^4\right), \quad S(t)-\frac{L_s}{2}<x<S(t)+\frac{L_s}{2},
\end{equation}
with $L_s = 0.5~m$ and $H_s = 0.026~m$ are being the slide height and thickness respectively. The movement of the slide is controlled by a function $S(t)$, given as:
\begin{equation}
    S(t) = \begin{cases}
        \frac{1}{2}a_0t^2 & ,0\leq t\leq t_1,\\
        \frac{1}{2}a_0t_1^2+u_t(t-t_1) & ,t_1< t\leq t_2,\\
        \frac{1}{2}a_0t_1^2+u_t(t-t_1)-\frac{1}{2}a_0(t-t_2)^2 & ,t_2< t\leq t_3,\\
        \frac{1}{2}a_0t_1^2+u_t(t_3-t_1)-\frac{1}{2}a_0(t_3-t_2)^2 & ,t>t_3.
    \end{cases}
\end{equation}
The comparison is made with three measured cases, following the parameters in Table \ref{tab:whittaker}.
\begin{table}[h!]
    \centering
    \begin{tabular}{|c|c|c|c|c|c|}
        \hline
         $Fr$ & $a_0~(m/s^2)$ & $u_t~(m/s)$ & $t_1~(s)$ & $t_2~(s)$ & $t_3~(s)$ \\
         \hline
         0.125 & $1.500$ & $0.163$ & $0.109$ & $2.109$ & $2.218$  \\
         0.25 & $1.500$ & $0.327$ & $0.218$ & $2.218$ & $2.436$ \\
         0.375 & $1.500$ & $0.491$ & $0.327$ & $2.327$ & $2.654$  \\
         \hline
    \end{tabular}
    \caption{Parameters of the experimental setup.}
    \label{tab:whittaker}
\end{table}

\begin{figure}[h!]
    \centering
    \includegraphics[width=1\linewidth]{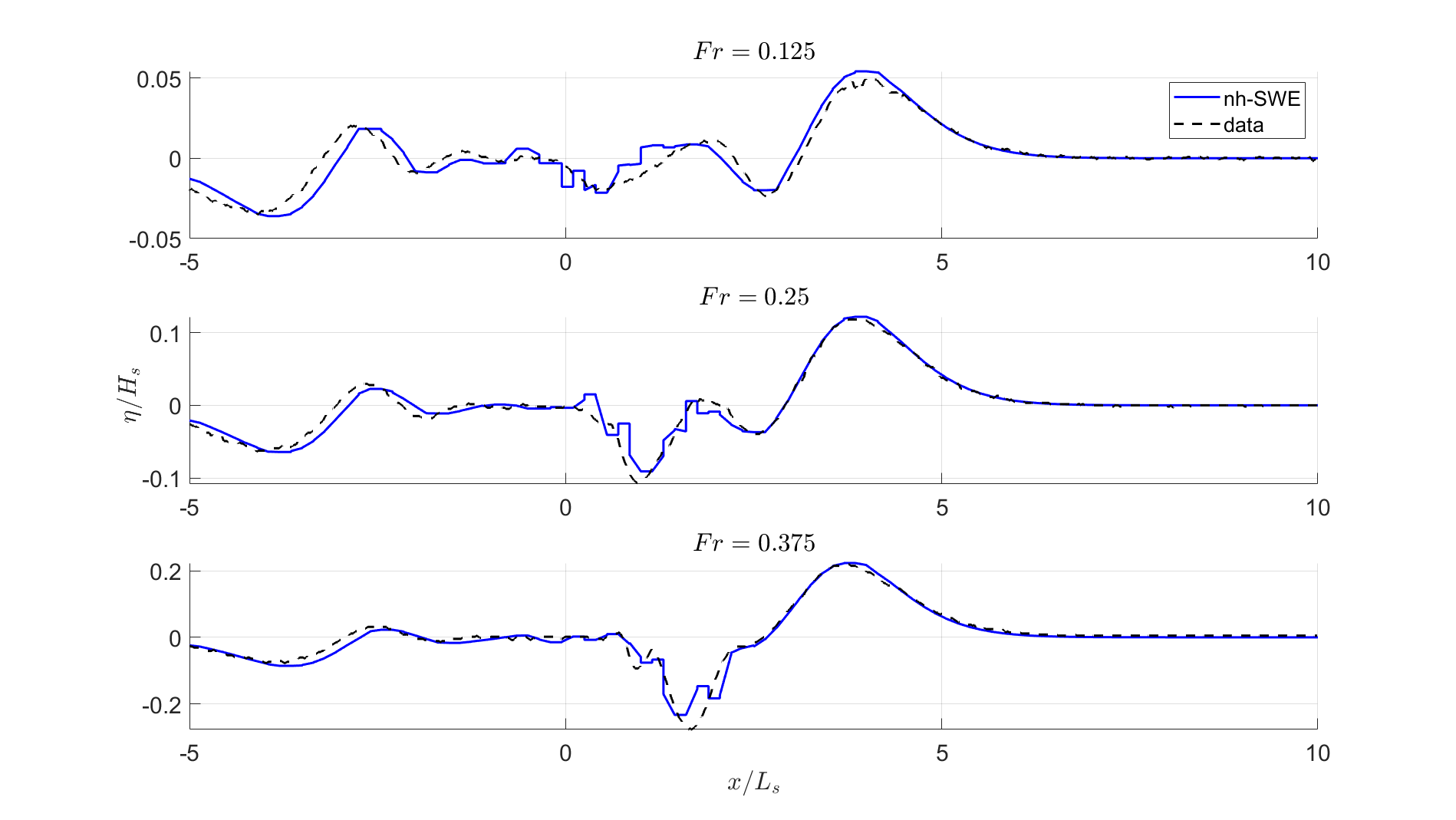}
    \caption{Comparison of the simulation result (blue solid line) with measured laboratory data (black dashed line) for a sliding semi-elliptic plate for $Fr=0.125$ (top), $0.25$ (middle), and $0.375$ (bottom).}
    \label{fig:global_whittaker}
\end{figure}

For our computational purpose, we discretize our model with $\Delta t = 0.005 ~s$ and $\Delta x = 0.075~m$. We compare our simulation with the measured data at $t = \frac{8}{\sqrt{g/L_s}}~s$ as shown in Fig. \ref{fig:global_whittaker}. Our simulations generally agree well with the data, which is represented by the produced errors and their correlation as in Table \ref{tab:global_whittaker}. 

\begin{table}[h!]
    \centering
    \begin{tabular}{|c|c|c|c|}
        \hline
         $Fr$ & RMSE & $r$ \\
         \hline
         0.125 & $0.0058$ & $0.9530$\\
         0.25  & $0.0095$ & $0.9750$\\
         0.375 & $0.0193$ & $0.9732$\\
         \hline
    \end{tabular}
    \caption{The produced error and correlation $r$ with the measured data for each experimental setting}
    \label{tab:global_whittaker}
\end{table}

\section{Towards Adaptive Model}\label{sec5}
\subsection{Study on Adaptivity Criteria}
In the previous part, we demonstrated the reliable and accurate modeling of weakly dispersive problems, including those generated by a moving bottom. As discussed earlier, the non-hydrostatic extension can be achieved by solving an elliptic system of equations based on the hydrostatic solution. Applying this correction globally results in an elliptic problem of the size equivalent to the number of all unknowns, requiring a correspondingly high computational cost. It is possible to solve this non-hydrostatic extension locally, thanks to the projection-based method. This means we restrict the non-hydrostatic region after the prediction step and only solve the elliptic equations in a sub-domain. To achieve that, a proper criterion is required to define this sub-domain. For assessing the accuracy, we will compare the adaptive model result similar to the previous evaluation. Moreover, we also measure the computational efficiency by comparing the computational time of the local and global application of the extension. 

For determining the area, where the non-hydrostatic extension is to be applied, it would be ideal to measure the size of this part of the equation, and apply it only where it exceeds a certain threshold. This estimate, however, needs to be derived from the hydrostatic predictor solution. One choice is to take $\frac{d}{\lambda}>0.05$ as a criterion, with $\lambda$ being the wavelength. We derive this from one of the main assumptions of shallow water theory, i.e., $\lambda \gg d$  \cite{dean1991}. However, determining the wavelength is not easy, especially for a more complicated case. Even a simple solitary wave case would have an infinite wavelength by definition. Hence, we take a first guess criterion based on the hydrostatic surface elevation $\tilde{\eta}$ given as $\left|\frac{\tilde{\eta}}{d}\right|$, another basic assumption of shallow water theory. Furthermore, we investigate any other possible criteria such as $|\tilde{\eta}_x|,\ |\tilde{u}|,\ |\tilde{u}_x|,\ |\tilde{w}|,$ and $|\tilde{w}_x|$, which all can be obtained from the predictor step. For further notation, we neglect the tilde for simplicity. With these criteria, we should only solve the elliptic problem on the area where the value of the respective criterion exceeds a threshold $k_{nh}$. We set $k_{nh} = 0.001$ for all our test cases.

\subsubsection{Solitary wave}
\begin{figure}[h]
    \centering
    \subfigure{%
        \includegraphics[width=0.4\textwidth]{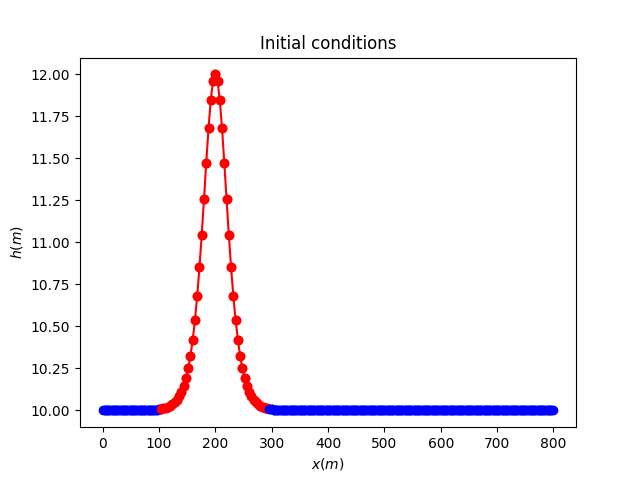}
    }
    \subfigure{%
        \includegraphics[width=0.4\textwidth]{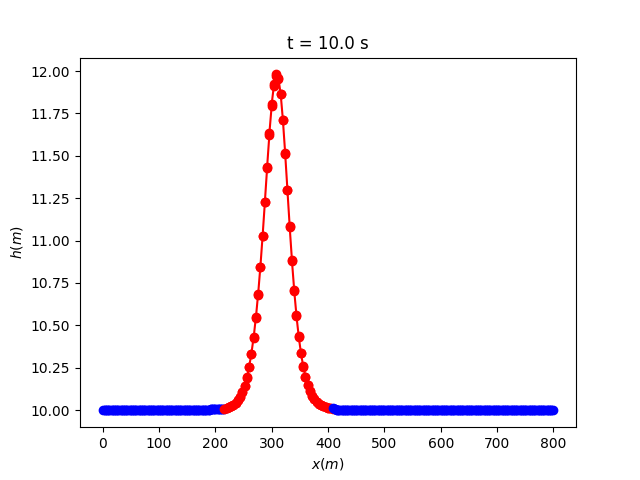}
    }
    \subfigure{%
        \includegraphics[width=0.4\textwidth]{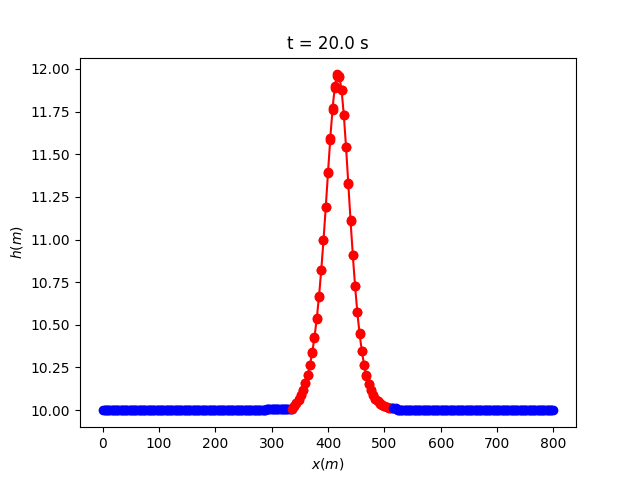}
    }
    \subfigure{%
        \includegraphics[width=0.4\textwidth]{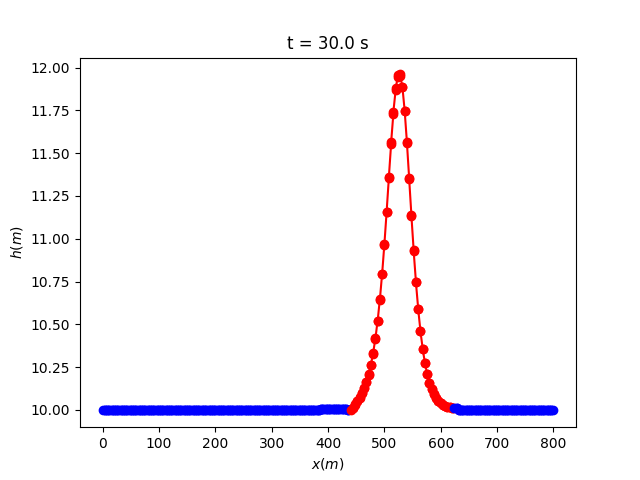}
    }
    \caption{Illustration of adaptive simulation of a propagating solitary wave, with splitting criterion $\left|\frac{\tilde{\eta}}{d}\right|>k_{nh} = 0.001$. The blue area depicts the hydrostatic area, while the red represents the non-hydrostatic domain, where we apply the correction.}
    \label{fig:local_solitary}
\end{figure}

\begin{figure}[h]
    \centering
    \subfigure{%
        \includegraphics[width=0.4\textwidth]{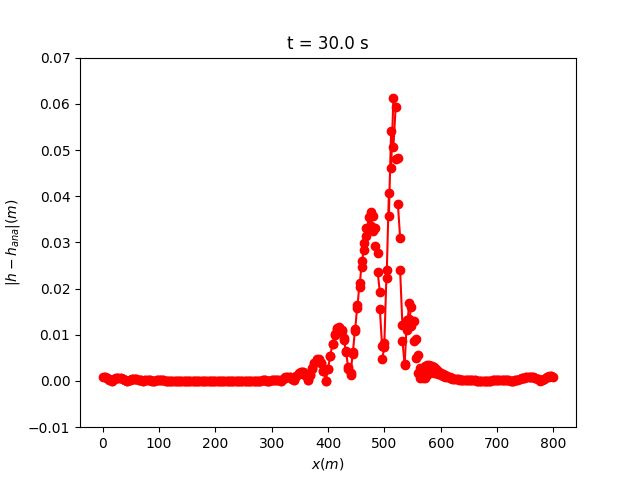}
    }
    \subfigure{%
        \includegraphics[width=0.4\textwidth]{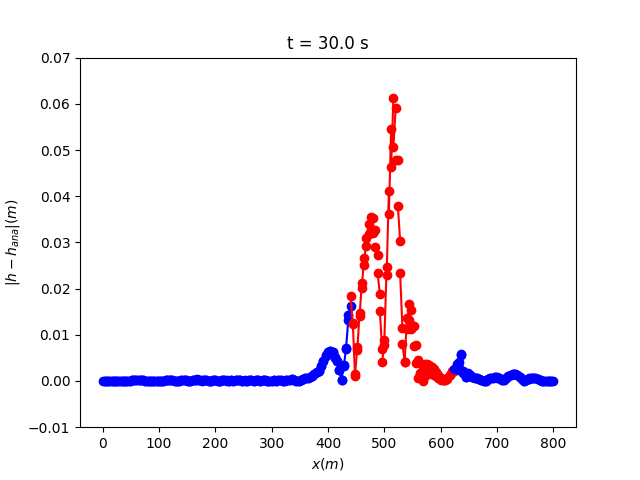}
    }
    \caption{Comparison of the absolute error of the global (left) and local (right) approach on the propagating solitary wave test case, with splitting criterion $\left|\frac{\tilde{\eta}}{d}\right|>k_{nh} = 0.001$.}
    \label{fig:error_solitary}
\end{figure}

To illustrate how our proposed adaptive model works, we begin with the simple solitary wave case as investigated in the former section.  Fig. \ref{fig:local_solitary} shows snapshots of an adaptive simulation on a propagating solitary wave with adaptive criterion $\left|\frac{\eta}{d}\right|>k_{nh}$. Physically, this criterion yields application of the correction to the area where significant elevation exists, which is intuitively reasonable for the solitary wave.
From this example, we can see how the proposed approach can reduce the computational effort with a smaller corrected region. Since an analytical solution is known, Fig. \ref{fig:error_solitary} shows the absolute errors of both global and local approaches at the end of the simulation time. It can be observed that the local adaptation manages to produce a similar accuracy. Differences are observable only at the interfaces between hydrostatic and non-hydrostatic regions. Note that the error is plotted for the discontinuous elements of the DG spatial discretization.

\begin{table}[h]
    \centering
    \begin{tabular}{|c|c|c|c|}
        \hline
         Criterion & $t_{\text{local}}/t_{\text{global}}$ & RMSE & $r$\\
         \hline
         $|\eta/d|$ & $0.1048$ & $0.0105$ & $0.9997$\\
         $|\eta_x|$ & $0.1016$ & $0.0136$ & $0.9995$\\
         $|u|$      & $0.2738$ & $0.0108$ & $0.9997$\\
         $|u_x|$    & $0.1285$ & $0.0158$ & $0.9993$\\
         $|w|$      & $0.1646$ & $0.0123$ & $0.9996$\\
         $|w_x|$    & $0.0818$ & $0.1146$ & $0.9606$\\
         \hline
    \end{tabular}
    \caption{Comparison of time ratio, errors, and correlations between the locally applied method and the analytical solution of propagating solitary wave with $k_{nh} = 0.001$ being the threshold.}
    \label{tab:local_solitary}
\end{table}


\begin{figure}[h]
    \centering
    \includegraphics[width=0.5\textwidth]{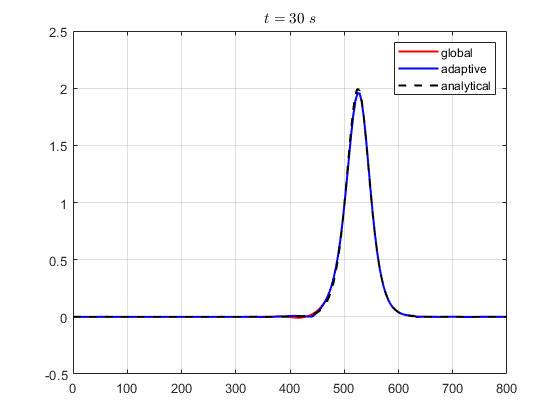}
    \caption{Comparison of the global simulation (red solid line), locally adaptive simulation (blue solid line) with the analytical solution (black dashed line), with criterion $\left|\frac{\tilde{\eta}}{d}\right|>k_{nh} = 0.001$}
    \label{fig:local_etad_solitary}
\end{figure}

To compare the performance of the proposed criterion, Table \ref{tab:local_solitary} lists the computational time ratio, the RMSE, and the correlation to the analytical solution at the end of the simulation time. It can be observed that most proposed criteria match the analytical solution, with $\left|\frac{\eta}{d}\right|$ giving the least error, while only taking $10.48\%$ of the global computational time (see Fig. \ref{fig:local_etad_solitary} for illustration). In contrast, criterion $|w_x|$ gives relatively poor results in both metrics (RMSE and $r$). This could be due to the fact that the form of the vertical velocity has two critical points (where $w_x = 0$), which leads to isolated elements. The existence of these isolated elements potentially introduces the larger error.

\subsubsection{Impulsive vertical thrust}
Note that in the next two cases, the initial vertical momentum value is zero, representing undisturbed water. However, this leads to static vertical momentum since the hydrostatic term also depends on the vertical momentum value from the previous step. Therefore, criteria $|w|$ and $|w_x|$ will not work. A slight modification is required for these two criteria, where we apply the model globally for one time-step and apply the model locally in the further steps. This approach allows a dynamic vertical momentum for the following steps.

\begin{figure}[h!]
    \centering
    \includegraphics[width=1\linewidth]{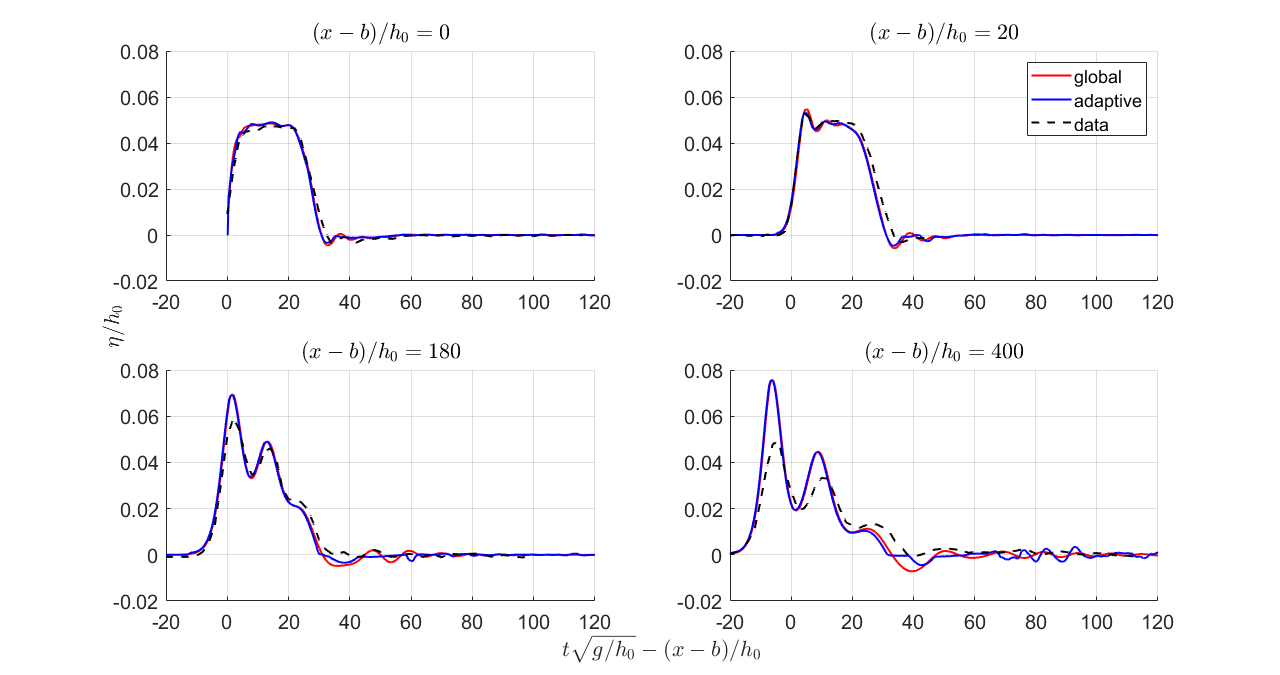}
    \caption{Comparison of the global simulation (red solid line), locally adaptive simulation (blue solid line), and measured laboratory data (black dashed line) recorded at $x = 0.61, 1.61, 9.61, 20.61~m$ for the uplift case.}
    \label{fig:local_up}
\end{figure}

\begin{figure}[h!]
    \centering
    \includegraphics[width=1\linewidth]{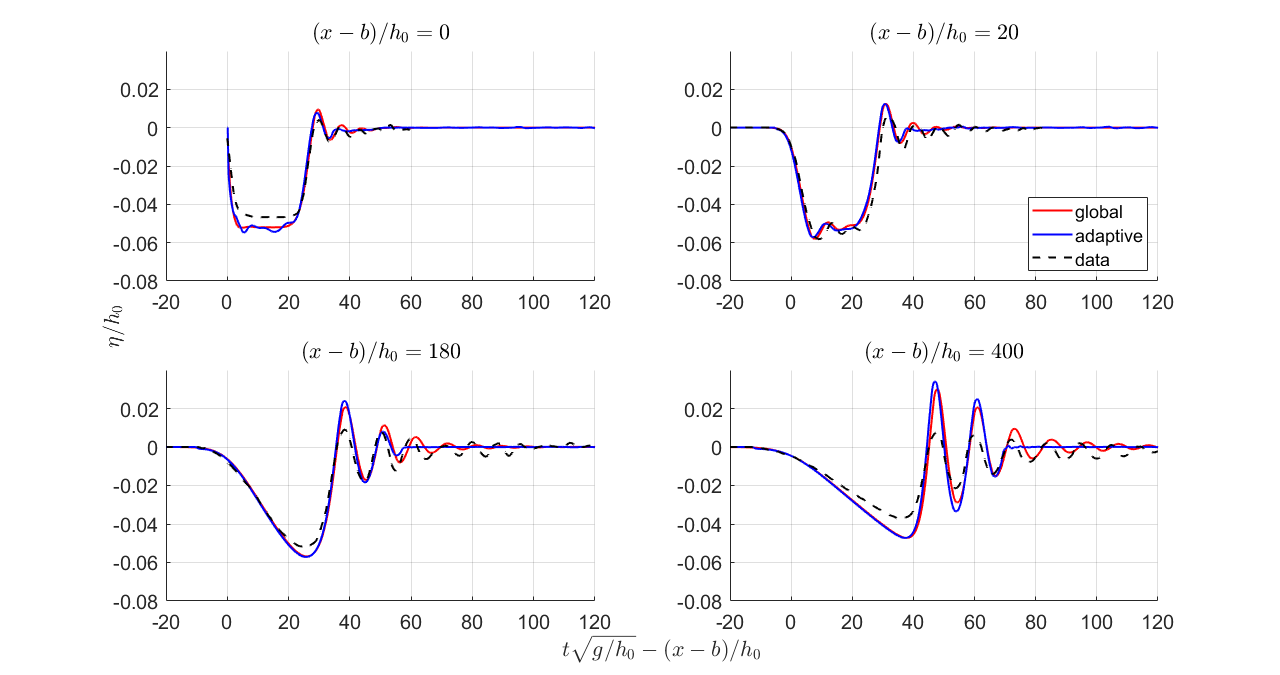}
    \caption{Comparison of the global simulation (red solid line), locally adaptive simulation (blue solid line), and measured laboratory data (black dashed line) recorded at $x = 0.61, 1.61, 9.61, 20.61~m$ for the downdraft case.}
    \label{fig:local_down}
\end{figure}

Table \ref{tab:local_upward} and \ref{tab:local_downward} show results from each criterion for upward and downward movement, respectively. By summing up the errors from each measurement point and considering their correlations, criterion $\left|\frac{\eta}{d}\right|$ gives the best results compared to others. Fig. \ref{fig:local_up} and \ref{fig:local_down} depict the results by taking such a criterion. The local approach with criterion $\left|\frac{\eta}{d}\right|$ manages to capture the significant dynamics while only taking $6.81\%$ and $8.79\%$ of the global computation time for the upward and downward movement, respectively.

\begin{table}[h!]
    \centering
    \resizebox{\textwidth}{!}{%
    \begin{tabular}{|c|c|c|c|c|c|c|c|c|c|}
        \hline
         \multirow{2}{*}{Criterion} & \multirow{2}{*}{$t_{\text{local}}/t_{\text{global}}$} & \multicolumn{4}{c|}{RMSE}& \multicolumn{4}{c|}{$r$}\\
         \cline{3-10}
         & & $\frac{x-b}{h_0} = 0$ & $\frac{x-b}{h_0} = 20$ & $\frac{x-b}{h_0} = 180$ & $\frac{x-b}{h_0} = 400$ & $\frac{x-b}{h_0} = 0$ & $\frac{x-b}{h_0} = 20$ & $\frac{x-b}{h_0} = 180$ & $\frac{x-b}{h_0} = 400$\\
         \hline
          $|\eta/d|$ & $0.0681$ & $0.0018$ & $0.0023$ & $0.0029$ & $0.0066$
                                & $0.9955$ & $0.9929$ & $0.9876$ & $0.9503$\\
          $|\eta_x|$ & $0.0587$ & $0.0023$ & $0.0027$ & $0.0040$ & $0.0059$
                                & $0.9925$ & $0.9903$ & $0.9707$ & $0.9330$\\
          $|u|$      & $0.1016$ & $0.0022$ & $0.0027$ & $0.0034$ & $0.0067$
                                & $0.9931$ & $0.9898$ & $0.9819$ & $0.9436$\\
          $|u_x|$    & $0.0869$ & $0.0021$ & $0.0025$ & $0.0036$ & $0.0058$
                                & $0.9938$ & $0.9917$ & $0.9758$ & $0.9380$\\
          $|w|$      & $0.0185$ & $0.0019$ & $0.0027$ & $0.0051$ & $0.0066$
                                & $0.9948$ & $0.9903$ & $0.9534$ & $0.9161$\\
          $|w_x|$    & $0.0220$ & $0.0030$ & $0.0031$ & $0.0052$ & $0.0066$
                                & $0.9871$ & $0.9863$ & $0.9529$ & $0.9179$\\
          \hline
    \end{tabular}%
    }
    \caption{Comparison of computational time ratio, errors, and correlations between the locally applied method and the measured laboratory data for the uplift vertical moving plate test case with $k_{nh} = 0.001$ being the threshold.}
    \label{tab:local_upward}
\end{table}

\begin{table}[h!]
    \centering
    \resizebox{\textwidth}{!}{%
    \begin{tabular}{|c|c|c|c|c|c|c|c|c|c|}
        \hline
         \multirow{2}{*}{Criterion} & \multirow{2}{*}{$t_{\text{local}}/t_{\text{global}}$} & \multicolumn{4}{c|}{RMSE}& \multicolumn{4}{c|}{$r$}\\
         \cline{3-10}
         & & $\frac{x-b}{h_0} = 0$ & $\frac{x-b}{h_0} = 20$ & $\frac{x-b}{h_0} = 180$ & $\frac{x-b}{h_0} = 400$ & $\frac{x-b}{h_0} = 0$ & $\frac{x-b}{h_0} = 20$ & $\frac{x-b}{h_0} = 180$ & $\frac{x-b}{h_0} = 400$\\
         \hline
          $|\eta/d|$ & $0.0879$ & $0.0034$ & $0.0036$ & $0.0037$ & $0.0068$                         
                                & $0.9928$ & $0.9834$ & $0.9841$ & $0.9555$\\
          $|\eta_x|$ & $0.0536$ & $0.0040$ & $0.0041$ & $0.0042$ & $0.0063$                         
                                & $0.9888$ & $0.9780$ & $0.9708$ & $0.9225$\\
          $|u|$      & $0.2267$ & $0.0037$ & $0.0039$ & $0.0044$ & $0.0066$                         
                                & $0.9912$ & $0.9807$ & $0.9755$ & $0.9447$\\
          $|u_x|$    & $0.1408$ & $0.0037$ & $0.0037$ & $0.0038$ & $0.0060$ 
                                & $0.9915$ & $0.9825$ & $0.9765$ & $0.9326$\\
          $|w|$      & $0.0329$ & $0.0033$ & $0.0044$ & $0.0050$ & $0.0071$ 
                                & $0.9927$ & $0.9750$ & $0.9597$ & $0.9103$\\
          $|w_x|$    & $0.0370$ & $0.0043$ & $0.0037$ & $0.0050$ & $0.0069$ 
                                & $0.9855$ & $0.9829$ & $0.9584$ & $0.9120$\\
          \hline
    \end{tabular}%
    }
    \caption{Comparison of computational time ratio, errors, and correlations between the locally applied method and the measured laboratory data for the downward vertical moving plate test case with $k_{nh} = 0.001$ being the threshold.}
    \label{tab:local_downward}
\end{table}


\subsubsection{Sliding bump over flat bottom}
\begin{figure}[h!]
    \centering
    \includegraphics[width=0.99\linewidth]{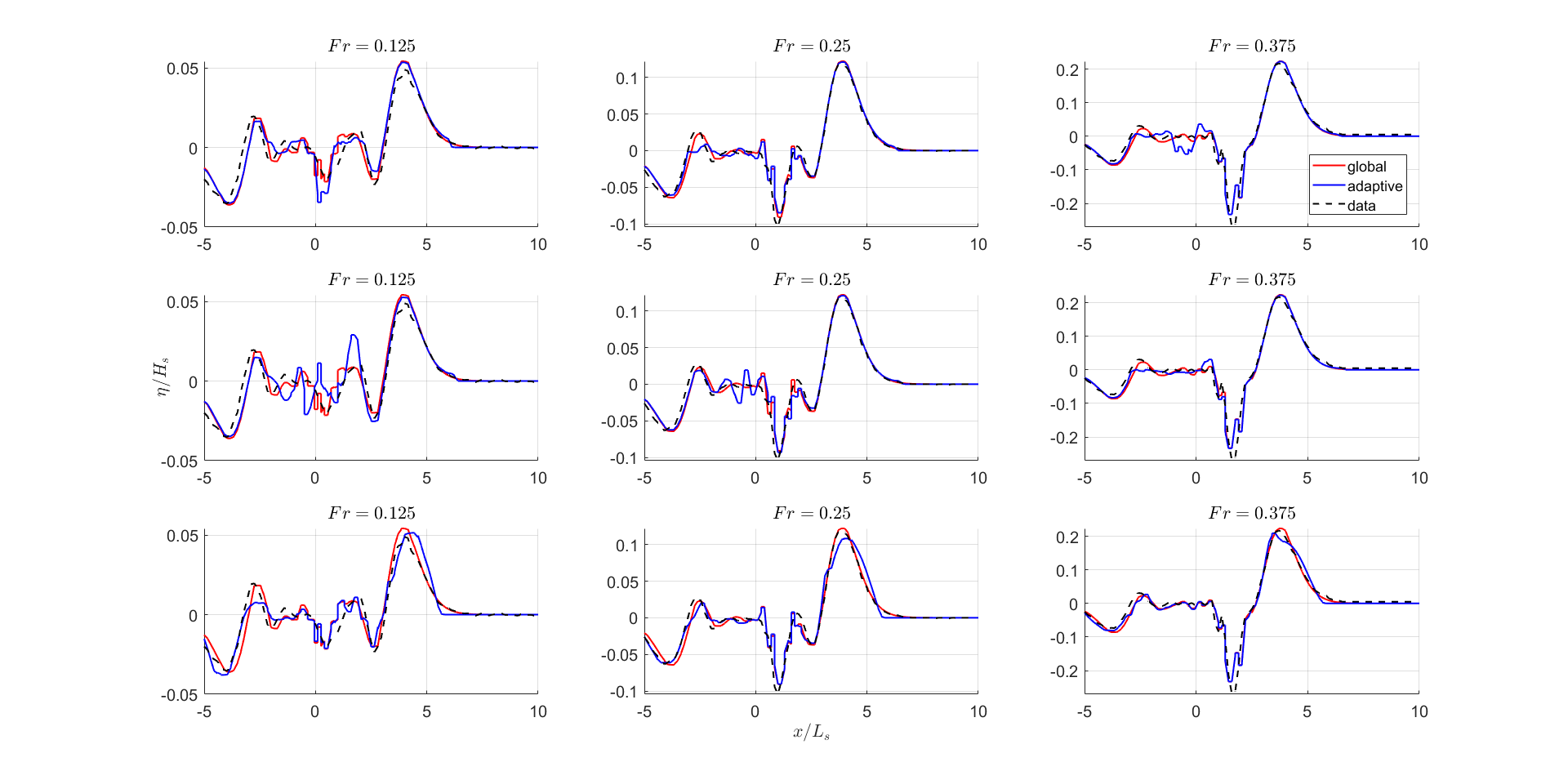}
    \caption{Comparisons of the global simulation (red solid line), locally adaptive model (blue solid line), with measured laboratory data (black dashed line) of the sliding bump test case, produced by criterion $\left|\frac{\eta}{d}\right|$ (top), $|u|$ (middle), and $|w_x|$ (bottom), with $k_{nh} = 0.001$ being the threshold.}
    \label{fig:local_whittaker}
\end{figure}

Finally, the adaptive simulation performance for the last test case can be seen in Table \ref{tab:local_whittaker}. In this case, however, the order of the performance of criteria $|w_x|$, $\left|\frac{\eta}{d}\right|$, and $|u|$ dominate these three cases interchangeably for both metrics. Out of curiosity, we illustrate comparing these three criteria in Fig. \ref{fig:local_whittaker}. It turns out that criterion $|w_x|$ mismatches the leading waveform, while the other two criteria match it well. However, it represents the following wave better, especially for a higher Froude number. In terms of computational time between these three criteria, criterion $|w_x|$ takes the shortest time, requiring $22.06-27.28\%$ of the global computational time. Following that, criterion $\left|\frac{\eta}{d}\right|$ takes $30.65-37.39\%$ of the global time, while $|u|$ requires $36.17-49.76\%$. 

\begin{table}[h]
    \centering
    \begin{tabular}{|c|c|c|c|}
        \multicolumn{4}{c}{$Fr = 0.125$}\\
        \hline
         Criterion & $t_{\text{local}}/t_{\text{global}}$ & RMSE & $r$\\
         \hline
         $|\eta/d|$ & $0.3065$ & $0.0064$ & $0.9435$\\
         $|\eta_x|$ & $0.2623$ & $0.0093$ & $0.8592$\\
         $|u|$      & $0.4976$ & $0.0079$ & $0.9130$\\
         $|u_x|$    & $0.4247$ & $0.0083$ & $0.8904$\\
         $|w|$      & $0.0965$ & $0.0205$ & $0.5696$\\
         $|w_x|$    & $0.2206$ & $0.0060$ & $0.9475$\\
         \hline
    \end{tabular}
    \hspace{1cm}
    \begin{tabular}{|c|c|c|c|}
        \multicolumn{4}{c}{$Fr = 0.25$}\\
        \hline
         Criterion & $t_{\text{local}}/t_{\text{global}}$ & RMSE & $r$\\
         \hline
         $|\eta/d|$ & $0.3624$ & $0.0103$ & $0.9703$\\
         $|\eta_x|$ & $0.2032$ & $0.0488$ & $0.6074$\\
         $|u|$      & $0.3977$ & $0.0108$ & $0.9675$\\
         $|u_x|$    & $0.3577$ & $0.0218$ & $0.8707$\\
         $|w|$      & $0.0862$ & $0.0290$ & $0.7885$\\
         $|w_x|$    & $0.2728$ & $0.0117$ & $0.9621$\\
         \hline
    \end{tabular}\\
    \vspace{0.5cm}
    \begin{tabular}{|c|c|c|c|}
        \multicolumn{4}{c}{$Fr = 0.375$}\\
        \hline
         Criterion & $t_{\text{local}}/t_{\text{global}}$ & RMSE & $r$\\
         \hline
         $|\eta/d|$ & $0.3739$ & $0.0224$ & $0.9635$\\
         $|\eta_x|$ & $0.2989$ & $0.0236$ & $0.9592$\\
         $|u|$      & $0.3617$ & $0.0201$ & $0.9711$\\
         $|u_x|$    & $0.3646$ & $0.0502$ & $0.8444$\\
         $|w|$      & $0.0967$ & $0.0473$ & $0.8441$\\
         $|w_x|$    & $0.2677$ & $0.0208$ & $0.9690$\\
         \hline
    \end{tabular}
    \caption{Comparison of computational time ration, errors, and correlations between the locally applied method and the measured laboratory data for sliding bump over flat bottom test case with $Fr = 0.125$ (top-left), $Fr = 0.25$ (top-right), and $Fr = 0.375$ (bottom), where $k_{nh} = 0.001$ being the threshold.}
    \label{tab:local_whittaker}
\end{table}

\subsection{Study on Improvement}
Motivated by the observation in the solitary test case, where the error in the adaptive run compared to the global run occurred at the interfaces, to improve the adaptive model further, we experiment by enlarging the correction domain by one element both to the left and right of the non-hydrostatic sub-domain. This allows to avoid the isolated elements issue, such as the one that clearly occurs in the solitary wave using criterion $|w_x|$. Moreover, using zero non-hydrostatic Dirichlet boundary conditions, where the interface is tight to a non-zero value of the non-hydrostatic extension, may lead to a higher error on the interface between the corrected and non-corrected regions. Therefore, the idea of extending the adaptivity area could relax this interaction. Below are the results of applying such an enlargement.

\subsubsection{Solitary wave}
Table \ref{tab:local_solitary(ex)} shows how enlarging the correction area by one element may overall improve the results by the proposed criteria, especially $|w_x|$. It can also be observed that the extension implies that all proposed criteria generally have similar accuracy to each other. Considering the computational time, criterion $\left|\frac{\eta}{d}\right|$ requires the least, taking only $15.93\%$ of the global computational time.
\begin{table}[h]
    \centering
    \begin{tabular}{|c|c|c|c|}
        \hline
         Criterion & $t_{\text{local}}/t_{\text{global}}$ & RMSE & $r$\\
         \hline
         $|\eta/d|$ & $0.1593$ & $0.0108$ & $0.9997$\\
         $|\eta_x|$ & $0.2555$ & $0.0108$ & $0.9997$\\
         $|u|$      & $0.3820$ & $0.0108$ & $0.9997$\\
         $|u_x|$    & $0.2946$ & $0.0118$ & $0.9996$\\
         $|w|$      & $0.3953$ & $0.0109$ & $0.9997$\\
         $|w_x|$    & $0.2884$ & $0.0110$ & $0.9997$\\
         \hline
    \end{tabular}
    \caption{Comparison of computational time ratio, errors, and correlations between the enlarged locally applied method with the analytical solution of propagating solitary wave with $k_{nh} = 0.001$ being the threshold.}
    \label{tab:local_solitary(ex)}
\end{table}

\subsubsection{Impulsive vertical thrust}
This test case again shows how extending the correction domain can improve the accuracy of the proposed adaptive model, especially on the produced waveform, as shown in Table \ref{tab:local_upward(ex)} and \ref{tab:local_downward(ex)} for upward and downward movement, respectively. Most proposed criteria give similar accuracies, with the exception of $|w|$, which clearly deviates based on the waveforms produced. This improvement comes with the cost of increased computational time, which is relatively significant compared to the previous test case due to the complexity of the dynamics. Criterion $\left|\frac{\eta}{d}\right|$ again demands the shortest computational time, with just $7.14\%$ and $11.63\%$ of the global application.
 \begin{table}[h!]
    \centering
    \resizebox{\textwidth}{!}{%
    \begin{tabular}{|c|c|c|c|c|c|c|c|c|c|}
        \hline
         \multirow{2}{*}{Criterion} & \multirow{2}{*}{$t_{\text{local}}/t_{\text{global}}$} & \multicolumn{4}{c|}{RMSE}& \multicolumn{4}{c|}{$r$}\\
         \cline{3-10}
         & & $\frac{x-b}{h_0} = 0$ & $\frac{x-b}{h_0} = 20$ & $\frac{x-b}{h_0} = 180$ & $\frac{x-b}{h_0} = 400$ & $\frac{x-b}{h_0} = 0$ & $\frac{x-b}{h_0} = 20$ & $\frac{x-b}{h_0} = 180$ & $\frac{x-b}{h_0} = 400$\\
         \hline
          $|\eta/d|$ & $0.0716$ & $0.0018$ & $0.0023$ & $0.0029$ & $0.0065$
                                & $0.9953$ & $0.9930$ & $0.9878$ & $0.9526$\\
          $|\eta_x|$ & $0.4117$ & $0.0021$ & $0.0025$ & $0.0031$ & $0.0067$
                                & $0.9939$ & $0.9916$ & $0.9850$ & $0.9433$\\
          $|u|$      & $0.4112$ & $0.0055$ & $0.0038$ & $0.0035$ & $0.0068$
                                & $0.9578$ & $0.9784$ & $0.9818$ & $0.9422$\\
          $|u_x|$    & $0.4397$ & $0.0019$ & $0.0024$ & $0.0030$ & $0.0066$
                                & $0.9947$ & $0.9922$ & $0.9869$ & $0.9505$\\
          $|w|$      & $0.3464$ & $0.0018$ & $0.0026$ & $0.0050$ & $0.0066$
                                & $0.9954$ & $0.9911$ & $0.9556$ & $0.9175$\\
          $|w_x|$    & $0.4299$ & $0.0022$ & $0.0025$ & $0.0028$ & $0.0063$
                                & $0.9930$ & $0.9916$ & $0.9894$ & $0.9605$\\
          \hline
    \end{tabular}%
    }
    \caption{Comparison of computational time ratio, errors, and correlations between the enlarged locally applied method and the measured laboratory data for the uplift vertical moving plate test case with $k_{nh} = 0.001$ being the threshold.}
    \label{tab:local_upward(ex)}
\end{table}

\begin{table}[h!]
    \centering
    \resizebox{\textwidth}{!}{%
    \begin{tabular}{|c|c|c|c|c|c|c|c|c|c|}
        \hline
         \multirow{2}{*}{Criterion} & \multirow{2}{*}{$t_{\text{local}}/t_{\text{global}}$} & \multicolumn{4}{c|}{RMSE}& \multicolumn{4}{c|}{$r$}\\
         \cline{3-10}
         & & $\frac{x-b}{h_0} = 0$ & $\frac{x-b}{h_0} = 20$ & $\frac{x-b}{h_0} = 180$ & $\frac{x-b}{h_0} = 400$ & $\frac{x-b}{h_0} = 0$ & $\frac{x-b}{h_0} = 20$ & $\frac{x-b}{h_0} = 180$ & $\frac{x-b}{h_0} = 400$\\
         \hline
          $|\eta/d|$ & $0.1163$ & $0.0034$ & $0.0036$ & $0.0037$ & $0.0068$
                                & $0.9932$ & $0.9832$ & $0.9848$ & $0.9602$\\
          $|\eta_x|$ & $0.6931$ & $0.0035$ & $0.0037$ & $0.0035$ & $0.0066$
                                & $0.9923$ & $0.9821$ & $0.9870$ & $0.9667$\\
          $|u|$      & $0.6204$ & $0.0035$ & $0.0038$ & $0.0038$ & $0.0069$
                                & $0.9921$ & $0.9817$ & $0.9830$ & $0.9545$\\
          $|u_x|$    & $0.8521$ & $0.0034$ & $0.0037$ & $0.0035$ & $0.0067$
                                & $0.9931$ & $0.9826$ & $0.9876$ & $0.9650$\\
          $|w|$      & $0.6139$ & $0.0031$ & $0.0033$ & $0.0055$ & $0.0069$
                                & $0.9948$ & $0.9866$ & $0.9516$ & $0.9105$\\
          $|w_x|$    & $0.6701$ & $0.0035$ & $0.0031$ & $0.0040$ & $0.0069$
                                & $0.9931$ & $0.9881$ & $0.9802$ & $0.9571$\\
          \hline
    \end{tabular}%
    }
    \caption{Comparison of computational time ratio, errors, and correlations between the enlargedlocally applied method and the measured laboratory data for the downdraft vertical moving plate test case with $k_{nh} = 0.001$ being the threshold.}
    \label{tab:local_downward(ex)}
\end{table}

\subsubsection{Sliding bump over flat bottom}
Finally, Table \ref{tab:local_whittaker(ex)} shows comparisons in our last test case of the sliding bottom. Expanding by one element homogenizes the results of most criteria with similar accuracy, with again the exception of $|w|$. The shortest simulation time was again achieved with criterion $\left|\frac{\eta}{d}\right|$ for all three different settings, which takes $37.6-40.32\%$ of the global computational time.
\begin{table}[h]
    \centering
    \begin{tabular}{|c|c|c|c|}
        \multicolumn{4}{c}{$Fr = 0.125$}\\
        \hline
         Criterion & $t_{\text{local}}/t_{\text{global}}$ & RMSE & $r$\\
         \hline
         $|\eta/d|$ & $0.3760$ & $0.0061$ & $0.9464$\\
         $|\eta_x|$ & $0.4135$ & $0.0057$ & $0.9552$\\
         $|u|$      & $0.5462$ & $0.0057$ & $0.9515$\\
         $|u_x|$    & $0.5855$ & $0.0059$ & $0.9500$\\
         $|w|$      & $0.2513$ & $0.0205$ & $0.5696$\\
         $|w_x|$    & $0.5445$ & $0.0058$ & $0.9531$\\
         \hline
    \end{tabular}
    \hspace{1cm}
    \begin{tabular}{|c|c|c|c|}
        \multicolumn{4}{c}{$Fr = 0.25$}\\
        \hline
         Criterion & $t_{\text{local}}/t_{\text{global}}$ & RMSE & $r$\\
         \hline
         $|\eta/d|$ & $0.4032$ & $0.0094$ & $0.9753$\\
         $|\eta_x|$ & $0.5118$ & $0.0095$ & $0.9751$\\
         $|u|$      & $0.4579$ & $0.0092$ & $0.9767$\\
         $|u_x|$    & $0.5973$ & $0.0096$ & $0.9743$\\
         $|w|$      & $0.2560$ & $0.0290$ & $0.7885$\\
         $|w_x|$    & $0.6147$ & $0.0095$ & $0.9749$\\
         \hline
    \end{tabular}\\
    \vspace{0.5cm}
    \begin{tabular}{|c|c|c|c|}
        \multicolumn{4}{c}{$Fr = 0.375$}\\
        \hline
         Criterion & $t_{\text{local}}/t_{\text{global}}$ & RMSE & $r$\\
         \hline
         $|\eta/d|$ & $0.4010$ & $0.0192$ & $0.9734$\\
         $|\eta_x|$ & $0.5005$ & $0.0193$ & $0.9733$\\
         $|u|$      & $0.5177$ & $0.0191$ & $0.9737$\\
         $|u_x|$    & $0.5948$ & $0.0191$ & $0.9738$\\
         $|w|$      & $0.2252$ & $0.0471$ & $0.8441$\\
         $|w_x|$    & $0.6176$ & $0.0193$ & $0.9732$\\
         \hline
    \end{tabular}
    \caption{Comparison of computational time ratio, errors, and correlations between the enlarged locally applied method and the measured laboratory data for sliding bump over flat bottom test case with $Fr = 0.125$ (top-left), $Fr = 0.25$ (top-right), and $Fr = 0.375$ (bottom), where $k_{nh} = 0.001$ being the threshold.}
    \label{tab:local_whittaker(ex)}
\end{table}

\section{Conclusions}\label{sec6}

This research shows that a non-hydrostatic adaptive model can be achieved with a simple criterion based on the solution of the hydrostatic SWE. Simply enlarging the adaptation area by one element leads us to more accurate and consistent results. The adaptive model allows us to attain a similar accuracy to the global model while reducing more than $50\%$ computational time. We plan to further investigate a more rigorous adaptive criterion and extend the applicability of this model to a two-dimensional and more realistic scenario. Involving the friction and viscosity terms are also of interest for further development.

\section*{Acknowledgements}
Part of this work is motivated by a Master’s thesis by Leila Wegener at Universit\"at Hamburg. The authors acknowledge the support of the Deutsche Forschungsgemeinschaft (DFG, German Research Foundation) within the Research Training Group GRK 2583 ``Modeling, Simulation and Optimization of Fluid Dynamic Applications''. J.B. additionally acknowledges funding by the DFG under Germany‘s Excellence Strategy – EXC 2037 ``CLICCS - Climate, Climatic Change, and Society'' – Project Number 390683824; as well as through the Collaborative Research Center TRR 181 ``Energy Transfers in Atmosphere and Ocean'' funded by the DFG - Project Number 274762653.

\section*{Conflicts of Interest}
The authors declare no conflicts of interest.

\section*{Data Availability Statement}
The data that support the findings of this study are available from the corresponding author upon reasonable request.

\bibliography{WileyNJD-AMA}
\end{document}